\documentclass[12pt,a4paper]{amsart}
\usepackage{graphicx}
\usepackage{epstopdf}
\usepackage{amssymb}
\usepackage{amsthm}
\usepackage{amsmath}
\usepackage{pstcol, pst-node}
\usepackage{mathrsfs}
\usepackage{amscd}
\usepackage[all]{xy}
\newtheorem{thm}{Theorem}

\newtheorem{lem}[thm]{Lemma}
\newtheorem{rem}[thm]{Remark}

\theoremstyle{definition}
\newtheorem{defn}[thm]{Definition}

\newtheorem{prop-def}[thm]{Proposition-Definition}

\theoremstyle{remark}

\begin{document}
\title{Rational curves in Fano hypersurfaces and tropical curves}
\author{Takeo Nishinou}
\date{}
\thanks{email : nishinou@math.tohoku.ac.jp}
\address{Mathematical Institute,
Tohoku University,
Sendai,
980-8578,
Japan } 
\subjclass[2000]{}
\maketitle
\begin{abstract}
Using ideas from the theory of tropical curves and degeneration, we
 prove that any Fano hypersurface (and more generally
 Fano complete intersections)
 is swept by at most quadratic rational curves.
\end{abstract}
\section{Introduction}
We begin with recalling a famous theorem of S.Mori \cite{Mo}.
\begin{thm}
Let $X$ be a nonsingular projective variety with $-K_X$ ample.
Then for any point $x\in X$, there is a rational curve through $x$.\qed
\end{thm}
Although more than thirty years have passed, 
 it seems that essentially there is no other proof of 
 this theorem other than Mori's original one, which uses reduction to positive characteristic.
In this paper, we attempt to prove this for Fano hypersurfaces
 in the realm of characteristic zero.
Namely, we prove the following.
We assume $n\geq 3$.

\begin{thm}\label{thm:hypersurface}
Let $X$ be an irreducible hypersurface of degree $d\leq n$
 in $\Bbb P^n$.
If $d = n$, then $X$ is swept by quadratic rational curves.
If $d<n$, then $X$ is swept by lines.
\end{thm}

This is already known (\cite{Ko}, Exercises 4.4 and 4.10) but we give a simple proof
 using tropical geometry, which may be of independent interest.
Our tools are ideas from the theory of tropical curves \cite{M, N, NS}
 and a calculation in the degenerate setting \cite{N2}.
Using degeneration, it suffices to consider rational curves in a projective space,
 and we further reduce the problem to a combinatorial one of finding
 a rational tropical curve with appropriate properties.

In the next section, we recall definitions and facts about tropical curves.
The proof of the theorem is given in Section \ref{sec:proof}.\\

\noindent
{\bf Acknowledgment.}
The author is supported by Grant-in-Aid for Young Scientists
(No.22740031).

\section{Tropical curves}
Here we quickly recall some definitions and facts about tropical curves.
See \cite{M, N, NS} for more details.

Let $\overline \Gamma$ be a weighted, connected finite graph.
Its sets of vertices and edges are denoted $\overline \Gamma^{[0]}$,
 $\overline \Gamma^{[1]}$, and $w_{\overline \Gamma} : 
  \overline \Gamma^{[1]} \to \Bbb N \setminus \{ 0 \}$
 is the weight function.
An edge $E \in \overline \Gamma^{[1]}$ has adjacent vertices
 $\partial E = \{ V_1, V_2 \}$.
Let $\overline \Gamma^{[0]}_{\infty} \subset \overline \Gamma^{[0]}$
 be the set of one-valent vertices.
We write $\Gamma = \overline \Gamma \setminus \overline\Gamma^{[0]}_{\infty}$.
Noncompact edges of $\Gamma$ are called \emph{unbounded edges}.
Let $\Gamma^{[1]}_{\infty}$ be the set of unbounded edges.
Let $\Gamma^{[0]}, \Gamma^{[1]}, w_{\Gamma}$
 be the sets of vertices and edges of $\Gamma$ and the weight function
 of $\Gamma$ (induced from $w_{\overline\Gamma}$ in an obvious way),
 respectively.
Let $N$ be a free abelian group of rank $n\geq 2$
 and $N_{\Bbb R} = N\otimes_{\Bbb Z}\Bbb R$.

\begin{defn}[{\cite[Definition 2.2]{M}}]\label{def:param-trop}
A \emph{parametrized tropical curve} in $N_{\Bbb R}$ is a proper map
 $h : \Gamma \to N_{\Bbb R}$ satisfying the following conditions.
\begin{enumerate}
\item[(i)] For every edge, $E \subset \Gamma$ the restriction $h \big|_E$
 is an embedding with the image $h(E)$ 
 contained in an affine line with rational slope, or $h(E)$ is a point.
\item[(ii)] For every vertex $V \in \Gamma^{[0]}$, $h(V)\in N_{\Bbb Q}$ and
 the following \emph{balancing
 condition} holds.
 Let $E_1, \dots, E_m \in \Gamma^{[1]}$ be the edges adjacent to $V$ and
 let $u_i \in N$ be the primitive integral vector emanating from $h(V)$
 in the direction of $h(E_i)$.
 Then
\begin{equation}
\sum_{j=1}^m w(E_j)u_j = 0.
\end{equation}
\end{enumerate}
\end{defn}

An isomorphism of parametrized tropical curves $h : \Gamma \to N_{\Bbb R}$ and 
 $h' : \Gamma' \to N_{\Bbb R}$ is a homeomorphism $\Phi : \Gamma \to \Gamma'$
 respecting the weights such that $h = h' \circ \Phi$.
\begin{defn}\label{def:tropical curve}
A \emph{tropical curve} is an isomorphism class of parametrized 
 tropical curves.
A tropical curve is \emph{trivalent} if $\Gamma$ is a trivalent graph.
The \emph{genus} of a tropical curve is the first Betti number of $\Gamma$.
The set of \emph{flags} of $\Gamma$ is 
\[
F\Gamma = \{(V, E) \big|
     V \in \partial E \}.
     \]
\end{defn}
\noindent
We often call a tropical curve of genus zero as a \emph{rational tropical curve}.
\begin{defn}\label{immersive}
We call a tropical curve $(\Gamma, h)$ \emph{immersive} if $h$ is an immersion
 and if $V\in\Gamma^{[0]}$, then $h^{-1}(h(V)) = \{V\}$.
\end{defn}
\noindent
In this paper, all the tropical curves we consider are trivalent and immersive.

By (i) of Definition \ref{def:param-trop}, we have a map
 $u : F\Gamma \to N$ sending a flag $(V, E)$
 to the primitive integral vector $u_{(V, E)} \in N$
 emanating from $h(V)$ in the direction of $h(E)$. 
\begin{defn}\label{def:type}
The (unmarked) \emph{combinatorial type} of a tropical curve $(\Gamma, h)$
 is the graph $\Gamma$ 
 together with the map $u : F\Gamma \to N$.
We write this by the pair $(\Gamma, u)$.
\end{defn}

\begin{defn}\label{def:degree}
The \emph{degree} of a type $(\Gamma, u)$
 is a function $\Delta: N \setminus \{ 0 \}
  \to \Bbb N$
 with finite support defined by 
\begin{equation*}
 \Delta(\Gamma, u)(v):= \sharp \{ (V, E) \in F\Gamma |
    E \in \Gamma^{[1]}_{\infty}, w(E)u_{(V, E)} = v \} 
\end{equation*}
Let $e=|\Delta| = \sum_{v\in N\setminus\{0\}}\Delta(v)$.
This is the same as the number of unbounded edges 
 of $\Gamma$ (not necessarily of $h(\Gamma)$).
\end{defn}

\subsection{Toric varieties associated to tropical curves and pre-log curves in them}
\begin{defn}\label{toric}
A toric variety $X$ defined by a fan $\Sigma$ 
 is called \emph{associated to a tropical curve $(\Gamma, h)$}
 if the set of the 
 rays of $\Sigma$ contains the set of the rays spanned by the vectors in $N$
 which are contained in the 
 support of the degree map $\Delta: N\setminus\{0\}\to \Bbb N$
 of $(\Gamma, h)$.

%If $ \mathfrak E$
% is an unbounded edge of $h(\Gamma)$, there is an obvious unique divisor of $X$
% corresponding to it.
%We write it as $D_{\mathfrak E}$ and call it the \emph{divisor associated to the edge 
% $\mathfrak E$}.
\end{defn}
Given a tropical curve $(\Gamma, h)$ in $N_{\Bbb R}$, 
 we can construct a polyhedral decomposition $\mathscr P$ of $N_{\Bbb R}$
 such that $h(\Gamma)$ is contained in the 1-skeleton of $\mathscr P$
 (\cite[Proposition 3.9]{NS}).
Given such $\mathscr P$, we construct a degenerating family $\mathfrak X\to \Bbb C$
 of
 a toric variety $X$ associated to $(\Gamma, h)$ (\cite[Section 3]{NS}).
We call such a family a \emph{degeneration of $X$ defined respecting $(\Gamma, h)$}.
Let $X_0$ be the central fiber.
It is a union $X_0 = \cup_{v\in\mathscr P^{[0]}}X_v$
 of toric varieties intersecting along toric strata.
Here $\mathscr P^{[0]}$ is the set of the vertices of $\mathscr P$.
\begin{defn}[{\cite[Definition 4.1]{NS}}]\label{torically transverse}
Let $X$ be a toric variety.
A holomorphic curve $C\subset X$ is \emph{torically transverse}
 if it is disjoint from all toric strata of codimension greater than one.
A stable map $\phi: C\to X$ is torically transverse if $\phi^{-1}(int X)\subset C$
 is dense and $\phi(C)\subset X$ is a torically transverse curve. 
Here $int X$ is the complement of the union of toric divisors.
\end{defn}
\begin{defn}\label{def:pre-log}
Let $C_0$ be a prestable curve.
A \emph{pre-log curve} on $X_0$ is a stable map $\varphi_0: C_0\to X_0$
 with the following properties.
\begin{enumerate}
\item[(i)] For any $v$, the restriction $C\times_{X_0}X_v\to X_v$
 is a torically transverse stable map.
\item[(ii)] Let $P\in C_0$ be a point which maps to the singular locus of $X_0$.
Then $C$ has a node at $P$, and $\varphi_0$ maps the two branches
 $(C_0', P), (C_0'', P)$ of $C_0$ at $P$ to different irreducible components 
 $X_{v'}, X_{v''}\subset X_0$.
Moreover, if $w'$ is the intersection index 
 of the restriction $(C_0', P)\to (X_{v'}, D')$ with the toric divisor
 $D'\subset X_{v'}$, 
 and $w''$ accordingly for $(C_0'', P)\to (X_{v''}, D'')$,
 then $w' = w''$.
\end{enumerate}
\end{defn}
Suppose we are given a 
 torically transverse rational curve $\varphi: \Bbb P^1\to X$ in a toric variety.
Then there is some rational tropical curve $(\Gamma, h)$ with the following 
 properties.
\begin{itemize}
\item The toric variety $X$ is associated to $(\Gamma, h)$.
\item Let $\mathfrak X\to \Bbb C$ be a degeneration of $X$ respecting $(\Gamma, h)$
 and $X_0 = \cup_{v\in\mathscr P^{[0]}}X_v$ the central fiber.
Then there is a family of prestable curves
$\mathfrak C\to \Bbb C$ whose generic fiber is $\Bbb P^1$, and
 a family of stable maps over $\Bbb C$
\[
\Phi: \mathfrak C\to \mathfrak X,
\] 
 such that the restriction to $1\in\Bbb C$ is $\varphi$.
\item The restriction to $0\in\Bbb C$, $\varphi_0: C_0\to X_0$ is 
 \emph{maximally degenerate}
 (see below).
\item The tropical curve $(\Gamma, h)$ is the dual intersection graph of this 
 maximally degenerate rational curve.
\end{itemize}
Here, a pre-log curve $\varphi_0: C_0\to X_0$ is maximally degenerate if
 for any $v\in \mathscr P^{[0]}$, the projection
 $\pi_v: C_0\times_{X_0}X_v\to X_v$ satisfies the following properties:
\begin{itemize}
\item Let $D_v$ be the union of toric divisors of $D_v$.
When $\dim X\geq 3$, then $\pi_v^{-1}(D_v)$ is at most three points, and
 the image of $\pi_v$ is contained in the closure of the orbit of 
 a one or two dimensional subtorus of the torus acting on $X$
 (note that this torus also acts on each component of $X_0$).  
\item When $\dim X = 2$, then the case where the image of $\pi_v$ 
 is the union of the
 closures of transversally intersecting orbits of one dimensional subtori
 is also allowed.
\end{itemize}
Conversely, given a rational tropical curve, we can construct a
 maximally degenerate rational curve in $X_0$, and we can lift it to a 
 smooth rational curve in a generic fiber of  $\mathfrak X\to \Bbb C$
 (when $\dim X = 2$, a nodal rational curve).
See \cite{NS} for more information about these results and definitions.
See also Remark \ref{rem:correspondence} below.

\section{Rational curves in Fano hypersurfaces}\label{sec:proof}
Here we give a proof of Theorem \ref{thm:hypersurface}.
It suffices to prove the claim for a hypersurface of degree $d$ defined by a 
 generic polynomial $f$.
Consider the degeneration
\[
z_0z_1\cdots z_{d-1}+tf = 0,
\]
 where $z_i$ are homogeneous coordinates of $\Bbb P^n$.
The central fiber $X_0$ is a union of $d$ $\Bbb P^{n-1}$s,
 intersecting along toric divisors.
Due to the assumption that 
 the degree $d$ is less than $n+1$,
 each component of $X_0$ has a divisor which is not contained
 in other components (we call it a \emph{free toric divisor}).
We mainly argue the case where each component of $X_0$ has
 just one free toric divisor (i.e., $d = n$), since the other cases are easier.
 
Singular locus $S$ of the total space $\mathfrak X$ of the degeneration is 
 given by the equations
\[
z_i = z_j = f = t = 0,\;\; i\neq j.
\]
Let
\[
X_0 = \cup_{i=1}^d \Bbb P_i^{n-1}
\]
 be the decomposition to irreducible components.
In $\Bbb P_i^{n-1}$, consider a rational curve of degree $d$.
If it is smooth, then it has $dn+n-4$ dimensional moduli.
As we argued in \cite{N}, a necessary condition for a curve 
\[
\varphi_0: C\to \Bbb P_i^{n-1}
\]
 to be liftable to a general fiber of $\mathfrak X\to \Bbb C$ is that
 any intersection of $\varphi_0(C)$ with the toric divisors of 
 $\Bbb P_i^{n-1}$ is contained in $S$.
This condition gives at most
 $d(n-1)$ dimensional condition.
We call them \emph{incidence condition}.
Here the factor of $n-1$ is the maximal number of non-free toric divisors.

Thus, the expected dimension of smooth rational curves of degree $d$
 satisfying the incidence condition is
\[
dn+n-4-d(n-1) = n+d-4.
\]
To sweep the hypersurface, we need at least $n-2$ dimensional family.
So we consider rational curves of degree two.
It suffices to prove that the obstruction cohomology class of a general member 
 of such a family vanishes.
 
Recall that a general embedded rational curve in $\Bbb P^{n-1}$
 can be described using a tropical curve (\cite{NS}).
We use this description to calculate the obstruction cohomology class,
 as we did in \cite{N, N2}.

By perturbing the incidence condition if necessary, we can assume
 the rational curve is generic, so that it corresponds to a trivalent,
 embedded tropical curve.
\begin{rem}\label{rem:correspondence}
By the statement that "rational curves in $\Bbb P^{n-1}$ are described by tropical curves", 
 we mean the following:
\begin{itemize}
\item Fix a general rational curve in $\Bbb P^{n-1}$ and a general rational
 tropical curve in $\Bbb R^{n-1}$ of the same degree (we do not impose any relation
 between these objects).
\item Then, we can take 
\begin{itemize}
\item a neighborhood $U$ of the rational curve in the moduli space of 
 rational curves of the given degree in $\Bbb P^{n-1}$,
 and 
\item a neighborhood $V$ of the rational tropical curve in the moduli space of 
 rational tropical curves of the given degree in $\Bbb R^{n-1}$
 ($V$ can be taken so that it is diffeomorphic to an open subset of $\Bbb R^N$
 for some $N$),
\end{itemize}
\item so that $U$ can naturally be considered as a complexification of $V$.
\end{itemize}
In particular, given a tropical curve, we cannot tell what the rational curve 
 precisely corresponding to it is.
In other words,
 there is no canonical correspondence between holomorphic and tropical curves
 (we need to specify some artificial incidence 
  conditions to obtain a precise correspondence between these two objects).
 But
 for our purpose, it suffices to know the following:
\begin{itemize}
\item Given a trivalent rational tropical curve of degree $d$,
 we can construct a degeneration 
 $\mathfrak P\to \Bbb C$ of $\Bbb P^{n-1}$ and also construct a
 $(dn+n-4)$-dimensional
 family $\mathcal F$ of
 maximally degenerate curves
 in the 
 central fiber $P_0$.
\item We can lift the curves in $\mathcal F$, so that we obtain
 a $(dn+n-4)$-dimensional
 family of smooth
 rational curves in a generic fiber of $\mathfrak P$.
%Note that $dn+n-4$ is the same as the dimension of the moduli space of 
% rational curves of degree $d$ in $\Bbb P^{n-1}$.
\item By suitably choosing $f$, some member of this $(dn+n-4)$-dimensional
 family satisfies the incidence condition imposed by $S$.
\end{itemize}
See \cite{NS} for the construction of the family of rational curves
 from tropical curves.
The last claim can be seen by a simple dimension count.
\end{rem}

\begin{figure}[h]
\includegraphics{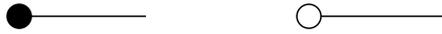}
\caption{Two types of unbounded edges with one-valent vertex attached.
The left figure represents the unbounded edge corresponding to the
 intersection with a non-free toric divisor, 
 and the right figure represents the unbounded edge corresponding to the
 intersection with a free toric divisor.}\label{fig:0}
\end{figure}

A usual tropical curve, as defined in Definition \ref{def:tropical curve}, is a proper map
\[
\varphi: \Gamma\to \Bbb R^{n-1}
\]
 from an abstract graph to the affine space.
In particular, it has unbounded edges, which correspond to the intersections with 
 the toric divisors.
However, in our situation, there are two types of
  intersections with the toric divisors: one free toric divisor and 
  $(n-1)$ non-free toric divisors.
So we add two types of one-valent vertices to unbounded edges to 
 distinguish them, see Figure \ref{fig:0}.

Now we describe the normal sheaf of a nodal rational curve corresponding to 
 a tropical curve.
We take a degeneration of $\Bbb P^{n-1}_i$ so that the rational curve
 becomes maximally degenerate (see \cite{NS}),
 and study the normal sheaf on each component.
The singular locus $S$ also degenerates, and we write its degeneration by the
 same letter $S$ (see \cite{N2}, Subsection 3.2.3).
A similar calculation was done in \cite{N}, the difference here is the existence of
 the non-free toric divisors.
 
A component of a maximally degenerate curve corresponds to a
 trivalent vertex of the corresponding tropical curve.
In particular, each component is contained in the closure of the orbit
 of a two dimensional subtorus of the torus acting on $\Bbb P_i^{n-1}$.
By perturbing the incidence condition again, we can assume the following:\\
 
 $(\ast)$ \emph{At each intersection of the 
 maximally degenerate rational curve and $S$, this orbit closure of 
 the two dimensional torus is transversal to $S$.}\\

There are three types of edges:
\begin{itemize}
\item An edge corresponding to a node.
\item An edge corresponding to an intersection with a free toric divisor.
\item An edge corresponding to an intersection with a non-free toric divisor.
\end{itemize}
If we do not take into account any special conditions
 (in other words, if each of the three edges is the
 one corresponding to an intersection with a free toric divisor),
 the normal sheaf of a
 torically transverse rational curve in an $n-1(\geq 2)$ dimensional toric variety
 which is contained in the orbit closure of a two dimensional subtorus
 is given by
\[
\mathcal O(1)\oplus \mathcal O^{\oplus n-3}.
\]
Here the component $\mathcal O(1)$ is the normal sheaf as a map to the orbit closure, 
 and the component $\mathcal O^{\oplus n-3}$ is the part transverse to the orbit closure.
By Serre duality, the first cohomology of it is dual to
\[
H^0(\Bbb P^1, (\mathcal O(-1)\oplus \mathcal O^{\oplus n-3})\otimes \omega_{\Bbb P^1}),
\]
 where $\omega_{\Bbb P^1}$ is the canonical sheaf.
This gives the (dual of) obstruction class, and we will calculate it 
 when there are above extra conditions.

The two from the three types of edges affect this calculation as follows:
\begin{itemize}
\item A node changes $\omega_{\Bbb P^1}$ to $\omega_{\Bbb P^1}(1)$
 (Serre duality for nodal curves).
\item An edge corresponding to an intersection with a non-free toric divisor
 changes $\mathcal O(1)$ component to $\mathcal O$, 
 by the calculation in \cite{N} and the above assumption $(\ast)$.
\end{itemize}
Our purpose was to show that there is a family of rational curves satisfying the 
 incidence condition $S$, such that the obstruction cohomology class of 
 a general member of it vanishes.
In view of Remark \ref{rem:correspondence}, it suffices to find a tropical curve whose
 obstruction cohomology class
 (more precisely, the obstruction cohomology class of
 the nodal rational curve associated to the tropical curve)
 calculated according to the above rule
 vanishes.

Such a tropical curve is given in the following way.
Let
\[
e_1 = (1, 0, \dots, 0),\;\; e_2 = (0, 1, 0, \dots, 0),\;\; \dots,\;\; e_{n-1} = (0, \dots, 0, 1)
\]
 be the standard basis of $\Bbb R^{n-1}$.
A rational tropical curve of degree $n-1$ in $\Bbb R^{n-1}$ 
 has $(n-1)$ unbounded edges in each directions
\[
-e_1, \;\; -e_2, \;\;\dots, \;\; -e_{n-1},\;\; d_n = e_1+e_2+\cdots e_{n-1},
\]
 where we take the direction of an unbounded edge to be the one emanating from 
 the unique adjacent vertex.
We assume $-e_1$ is the direction corresponding to the (unique) free face.

Then consider the following rational quadratic tropical curve.
Here and hereafter, all the edges have weight one.

\begin{figure}[h]
\includegraphics{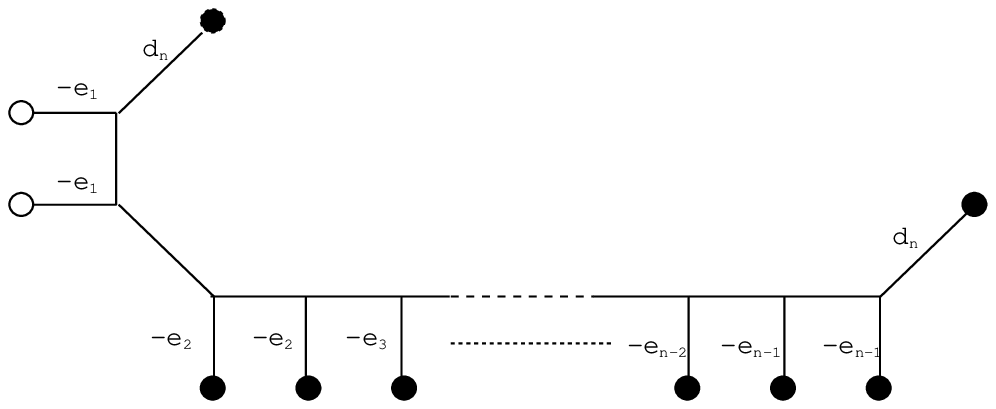}
\caption{}\label{figure:1}
\end{figure}

Now we prove the following, which completes the proof of Theorem \ref{thm:hypersurface}:
\begin{lem}
The maximally degenerate
 curve corresponding to the tropical curve of Figure \ref{figure:1}
 has vanishing obstruction.
\end{lem}
\proof
The tropical curve has four types of trivalent vertices:

\begin{figure}[h]
\includegraphics{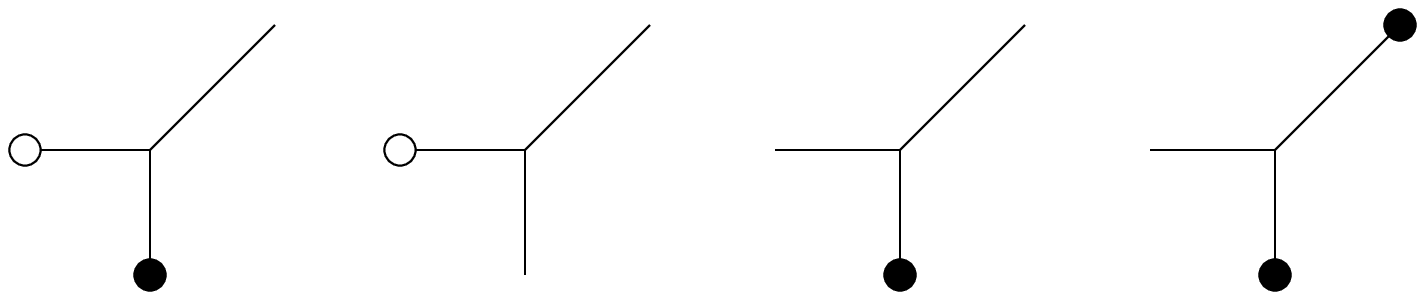}
\caption{}\label{}
\end{figure}

According to the rule described above, 
 the normal sheaf of the component corresponding to the leftmost vertex is given by
\[
\mathcal O\oplus \mathcal O^{n-3}.
\]
By Serre duality for nodal curves, the obstruction cohomology of this component is
 isomorphic to the dual of the following:
\[
H^0(\Bbb P^1, (\mathcal O\oplus\mathcal O^{n-3})\otimes \omega_{\Bbb P^1}(1))
 = H^0(\Bbb P^1, \mathcal O(-1)\oplus\mathcal O^{n-3}(-1))
 = 0.
\]

Similarly, the normal sheaves  are
 given by the following:
\[
\mathcal O(1)\oplus \mathcal O^{n-3},
\]
\[
\mathcal O\oplus\mathcal O^{n-3},
\]
\[
\mathcal O(-1)\oplus\mathcal O^{n-3}.
\]
Their corresponding (dual of the) obstruction classes are given as follows:
\[
H^0(\Bbb P^1, (\mathcal O(-1)\oplus\mathcal O^{n-3})\otimes \omega_{\Bbb P^1}(2))
 = H^0(\Bbb P^1, \mathcal O(-1)\oplus\mathcal O^{n-3}),
\]
\[
H^0(\Bbb P^1, (\mathcal O\oplus\mathcal O^{n-3})\otimes \omega_{\Bbb P^1}(2))
 = H^0(\Bbb P^1, \mathcal O\oplus\mathcal O^{n-3}),
\]
\[
H^0(\Bbb P^1, (\mathcal O(1)\oplus\mathcal O^{n-3})\otimes \omega_{\Bbb P^1}(1))
 = H^0(\Bbb P^1, \mathcal O\oplus\mathcal O^{n-3}(-1)).
\]
In each of these three cases, the cohomology does not vanish, but one sees that 
 if one knows that a section representing the cohomology class is zero at one point, 
 then that section is itself zero.

Now let us look at the tropical curve of Figure \ref{figure:1}.
There is a vertex of the leftmost type,
 and the obstruction class restricted to 
 this component is zero.
Then by above observation, one sees that the obstruction of 
 the whole maximally degenerate curve corresponding to the tropical curve is 
 zero.
 
To see that the family of these curves sweeps the hypersurface, one has to check that
 the family is not contained in a proper subvariety.
However, it is easy to see that in the tropical curve of Figure \ref{figure:1},
 with the normal sheaf given above, one of the two vertices intersecting the free face
 has $n-2$ dimensional freedom to move, so sweeps an open subset of the 
 free face (once the place of this vertex is chosen, 
 the place of the other vertex is determined uniquely by incidence conditions).
This shows that the family of curves cannot contained in a proper subvariety. \qed

\begin{rem}
By the same proof as above and the calculation at the beginning of this section, 
 it is easy to see that a Fano hypersurface of degree $d\leq n-2$ in $\Bbb P^{n}$
 is swept by lines, which may be easily proved by induction taking 
 various hyperplane sections. 
\end{rem}
\begin{rem}
As in \cite{Katz}, the argument for the degeneration of hypersurfaces is 
 extended to complete intersection without any essential change.
Thus, a Fano complete intersection 
\[
X = \cap_{i=1}^k V_i\subset \Bbb P^{k+n-1},
\]
 where $V_i$ is a hypersurface of degree $n_i>1$, with
\[
\sum_{i = 1}^k n_i = d\leq k+n-1
\]
 is swept by quadratic rational curves (when $d = k+n-1$),
 or by lines (when $d<k+n-1$).
\end{rem}

\begin{rem}
Combining with the results in \cite{N} and moderate tropical intersection theory, 
 we can study more general higher genus curves in Fano 
 (or more general varieties which have toric degeneration) varieties
 using tropical technique.
\end{rem}


\begin{thebibliography}{99} 
\bibitem{Katz}{\sc Katz,S.,}
{\it Lines on complete intersection threefolds with $K=0$.}
Math. Z. 191 (1986), no. 2, 293-296.
\bibitem{Ko}{\sc Koll\'ar,J.,}
{\it Rational curves on algebraic varieties.}
Ergebnisse der Mathematik und ihrer Grenzgebiete. 3. 
 Folge A Series of Modern Surveys in Mathematics, Vol. 32.
\bibitem{M}{\sc Mikhalkin,G.,}
{\it Enumerative tropical algebraic geometry in $\Bbb R\sp 2$.}
  J. Amer. Math. Soc.  18  (2005),  no. 2, 313--377.


\bibitem{Mo}{\sc Mori,S.,}
{\it Projective manifolds with ample tangent bundles.} Ann. Math. 110 (3): 593-606.

\bibitem{N}{\sc Nishinou,T.,}
{\it Correspondence theorems for tropical curves.}
Preprint.
\bibitem{N2}{\sc Nishinou,T.,}
{\it Counting curves via degeneration.}
Preprint.


\bibitem{NS}{\sc Nishinou,T. and Siebert,B.,}
{\it Toric degenerations of toric varieties and tropical curves.}
 Duke Math. J. 135 (2006), no. 1, 1--51.
 
 
 
 
 
 
 
 
 
 
 
 
 
 
 
 
 
 
 


\end{thebibliography}
\end{document}